\documentclass[reqno,12pt]{amsart}

\newtheorem{Theorem}{Theorem}[section]

\newtheorem{Corollary}[Theorem]{Corollary}
\newtheorem{Lemma}[Theorem]{Lemma}
\theoremstyle{remark}
\newtheorem{Remark}[Theorem]{Remark}
\numberwithin{equation}{section}

\def\ta{\theta}

\begin{document}

\title[Elliptic matrix inversion and elliptic hypergeometric series]
{On Warnaar's elliptic matrix inversion \\and Karlsson--Minton-type\\
 elliptic  hypergeometric series}

\author{Hjalmar Rosengren}
\address{Department of Mathematics, Chalmers University of Technology and
G\"oteborg University, SE-412 96 G\"oteborg, Sweden}
\email{hjalmar@math.chalmers.se}
\urladdr{http://www.math.chalmers.se/{\textasciitilde}hjalmar}

\author[Michael Schlosser]{Michael Schlosser$^*$}
\address{Institut f\"ur Mathematik der Universit\"at Wien,
Strudlhofgasse 4, A-1090 Wien, Austria}
\email{schlosse@ap.univie.ac.at}
\urladdr{http://www.mat.univie.ac.at/{\textasciitilde}schlosse}

\thanks{$^*$The second author was supported by an APART grant of the Austrian
Academy of Sciences}
%\date{September 17, 2003}

\begin{abstract}
Using Krattenthaler's operator method,
we give a new proof of Warnaar's recent elliptic extension of
Krattenthaler's matrix inversion. Further, using a theta function
identity closely related to Warnaar's inversion,
 we derive  summation and transformation formulas for elliptic
hypergeometric series of Karlsson--Minton-type. A special case yields a
particular  summation that was used by Warnaar
to derive quadratic, cubic and quartic transformations for elliptic
hypergeometric series. 
Starting from another theta function identity, we derive yet different
summation and transformation formulas for elliptic hypergeometric series
of Karlsson--Minton-type. These latter identities seem quite unusual
and appear to be new already in the trigonometric (i.e., $p=0$) case.
\end{abstract}

\maketitle

\section{Introduction}\label{sec0}

Matrix inversions provide a fundamental tool for studying hypergeo\-metric
and basic hypergeometric (or $q$-) series. For instance, they underlie
the celebrated Bailey transform \cite{andrews}. For multiple
hypergeometric series, multidimensional matrix inversions have
similarly proved a powerful tool, see \cite
{bmil,chum,kschl,lasschl,lilmil,milne,millil,schlossmmi,schlnammi,schlnmmi}.

Recently, a new class of generalized hypergeometric series
 was introduced, the
elliptic  hypergeometric series of Frenkel and Turaev
\cite{frtu}. In \cite{warnell}, Warnaar
found an elliptic matrix inversion and used it to obtain several new
quadratic, cubic and quartic summation and transformation formulas for
elliptic hypergeometric series.

Warnaar's matrix inversion can
 be stated as follows \cite[Lemma~3.2]{warnell}.
If
\begin{subequations}
\label{fg}
\begin{equation}\label{fnk}
f_{nk}=\frac{\prod_{j=k}^{n-1}\ta(a_jc_k)\ta(a_j/c_k)}
{\prod_{j=k+1}^{n}\ta(c_jc_k)\ta(c_j/c_k)}
\end{equation}
and
\begin{equation}\label{gkl}
g_{kl}=
\frac{c_l\,\ta(a_lc_l)\ta(a_l/c_l)}{c_k\,\ta(a_kc_k)\ta(a_k/c_k)}
\frac{\prod_{j=l+1}^{k}\ta(a_jc_k)\ta(a_j/c_k)}
{\prod_{j=l}^{k-1}\ta(c_jc_k)\ta(c_j/c_k)},
\end{equation}
\end{subequations}
then the infinite lower-triangular matrices $(f_{nk})_{n,k\in\mathbb Z}$
and $(g_{kl})_{k,l\in\mathbb Z}$ are {\em inverses} of each other, i.~e.,
the orthogonality relations
\begin{equation}\label{wmi}
\sum_{k=l}^nf_{nk}g_{kl}=\delta_{nl}, \qquad\mbox{for all $n,l\in\mathbb Z$}
\end{equation}
and (equivalently)
\begin{equation}\label{pmi}
\sum_{k=l}^ng_{nk}f_{kl}=\delta_{nl}, \qquad\mbox{for all $n,l\in\mathbb Z$}
\end{equation}
hold. 
In \eqref{fnk} and \eqref{gkl}, $\ta(x)$ is the {\em theta function},
defined by
\begin{equation*}
\ta(x)=\ta(x;p):=\prod_{j=0}^\infty(1-xp^j)(1-p^{j+1}/x),
\end{equation*}
for $|p|<1$.

Note that $\ta(x)$ reduces to $1-x$ for $p=0$. 
In this case Warnaar's matrix inversion reduces to a result of
Krattenthaler \cite[Corollary]{krattmi},
which in turn generalizes a large number of previously known
explicit matrix inversions. 

The present paper can be viewed as a spin-off of an attempt to
obtain multivariable extensions of Warnaar's matrix inversion and use
these to study elliptic hypergeometric series related to classical
root systems. This  led us to discover several aspects of
Warnaar's result which are interesting already in the one-variable
case. Multivariable extensions of these ideas are postponed to
future publications. 

Warnaar's proof of his inversion is based on the equation \eqref{pmi},
which is obtained as a special case of a more general identity, the
latter being easily proved by induction. This approach seems 
difficult (though interesting) to generalize to the multivariable case. 
On the other hand,  as was pointed out in \cite{rosers},
 the  identity \eqref{wmi} for Warnaar's inversion is  equivalent to a
 partial fraction-type expansion for theta functions due to Gustafson,
 \eqref{dpf} below. 
This leads to a  short proof of Warnaar's
(and thus also Krattenthaler's) matrix inversion, which is
described in Section \ref{sec:pf}. 

In another direction, Krattenthaler's 
proof of the case $p=0$ 
used a certain ``operator method'', cf.\
 Lemma \ref{lemma} below.
In Section \ref{sec:op} we extend Krattenthaler's proof to the
elliptic case. This requires some non-obvious steps, essentially
because addition formulas for theta functions are more complicated
than those for trigonometric functions implicitly used by Krattenthaler.

 We hope that both  the elementary proof of Warnaar's inversion
 given in Section~\ref{sec:pf}
 and the operator proof given in Section  \ref{sec:op} will
be useful for finding multivariable extensions.

Apart from the  matrix inversion \eqref{fg}, another important tool in
Warnaar's paper is the identity \cite[Theorem 4.1]{warnell}
(see \eqref{esf} below for the notation), which we write as
\begin{multline}\label{wbb}
\sum_{k=0}^N\frac{\theta(aq^{2ks})}{\theta(a)}
\frac{(a,q^{-Ns},b,a/b;q^s)_k}{(q^s,aq^{(N+1)s},aq^s/b,bq^s;q^s)_k}
\frac{(cq^N,aq/c;q)_{sk}}
{(aq^{1-N}/c,c;q)_{sk}}\,q^{sk}
\\
=\frac{(aq^s,q^s;q^s)_N}{(bq^s,aq^s/b;q^s)_N}
\frac{(c/b,bc/a;q)_{N}}{(c,c/a;q)_{N}}.
 \end{multline}
Here, $s$ is a positive
and $N$ a non-negative integer. 
In \cite{warnell}, this
 was obtained by combining \eqref{wmi} for the inverse pair \eqref{fg}
with a certain bibasic summation.
 The identity \eqref{wbb} was then applied, with $s=2$, $3$ and $4$, to
obtain  quadratic, cubic and quartic elliptic hypergeometric
identities, respectively. 

A characteristic property of \eqref{wbb} is that certain quotients of
numerator and denominator parameters (such as $b/q^sb$  and $cq^N/c$) are
integral powers of $q$. Classical and basic hypergeometric series with 
the analogous property have been called
\emph{Karlsson--Minton-type} and ($q$-)IPD-type (for Integral Parameter
Differences)   series. A seminal result for such series
is Minton's summation formula \cite{minton}
\begin{equation}\label{msf}
{}_{r+2}F_{r+1}\left(\begin{matrix}-N,b,c_1+m_1,\dots,c_r+m_r\\
b+1,c_1,\dots,c_r\end{matrix}\,;1\right)=\frac{N!}{(b+1)_N}
\prod_{i=1}^r\frac{(c_i-b)_{m_i}}{(c_i)_{m_i}}, \end{equation}
where it is assumed that $m_i$ are non-negative integers with
$|m|:=\sum_im_i\leq N$.
This has been extended to non-terminating, bilateral and well-poised series  
\cite{chu0,chu,gasper1,gasper2,karlsson,schlossed} and further to multiple
series \cite{ros1,ros2,schlossmt}.
However, for elliptic hypergeometric series, \eqref{wbb} has until now
been an isolated result.

At first sight, \eqref{wbb} looks somewhat different from known
Karlsson--Minton-type identities. However, writing
\begin{equation}\label{xsk}(x;q)_{sk}=(x,xq,\dots,xq^{s-1};q^{s})_k, 
\end{equation}
 it is not hard to check that it can be obtained as a special case of
 the more conventional summation formula
\begin{multline}\label{kmsi}\sum_{k=0}^N\frac{\theta(aq^{2k})}{\theta(a)}
\frac{(a,q^{-N},b,a/b;q)_k}{(q,aq^{N+1},aq/b,bq;q)_k}\,q^k
\prod_{j=1}^r\frac{(c_jq^{m_j},aq/c_j;q)_k}{(aq^{1-m_j}/c_j,c_j;q)_k}\\
=\frac{(aq,q;q)_N}{(bq,aq/b;q)_N}\prod_{j=1}^r\frac{(c_j/b,c_jb/a;q)_{m_j}}
{(c_j,c_j/a;q)_{m_j}},\qquad |m|=N\end{multline}
(with $q$ replaced by $q^s$).  This result will be
  proved  in  Section \ref{sec:km}. When
$p=0$, \eqref{kmsi} reduces to a special case
of an identity of Gasper \cite[Eq.~(5.13)]{gasper2}, which in turn
contains \eqref{msf} as a degenerate case.

Gasper's proof of \eqref{kmsi} in the case  $p=0$ does not
immediately extend to the elliptic case. A different proof was given by Chu
\cite{chu}, who independently obtained and generalized Gasper's identity
by recognizing it as a special case of a partial fraction expansion.
In Section \ref{sec:km}  we use Chu's method to generalize 
\eqref{kmsi} in a different direction, namely, to 
 a  multiterm Karlsson--Minton-type
transformation, Theorem \ref{kmt}. It is obtained as a special case of
 Gustafson's  identity \eqref{dpf}, or equivalently of 
\eqref{wmi} for Warnaar's inversion.
 Theorem \ref{kmt} may be viewed as an elliptic
analogue of Sears' transformation for well-poised series, cf.\ Remark
\ref{searsremark}.

In Section \ref{sec:ekm}, we repeat the analysis of Section 4, starting from a
different elliptic partial fraction identity, \eqref{apf}. 
This leads to some  exotic summation and transformation 
formulas for Karlsson--Minton-type elliptic hyper\-geo\-metric series,
which appear to be new also when $p=0$. 

Finally, in the Appendix we give an alternative proof of \eqref{kmsi},
using induction on $N$. We hope that the two proofs we 
give of this identity will both be useful for finding multivariable 
extensions of \eqref{kmsi}, and of related quadratic, cubic and
quartic identities from \cite{warnell}.

{\bf Notation:}
We have already introduced the theta function
$\theta(x)=\theta(x;p)$. The nome $p$ is fixed throughout and will be
suppressed from the notation.
We sometimes write
\begin{equation}
\theta(x_1,\dots,x_n):=\theta(x_1)\dots\theta(x_n)
\end{equation}
for brevity. We will frequently use the following
two properties of theta functions:
\begin{equation}\label{tif}
\ta(x)=-x\,\ta(1/x)
\end{equation}
and the {\em addition formula}
\begin{equation}\label{tfaddf}
\ta(xy,x/y,uv,u/v)-\ta(xv,x/v,uy,u/y)
=\frac uy\,\ta(yv,y/v,xu,x/u)
\end{equation}
(cf.\ \cite[p.~451, Example 5]{whitwatson}).

We denote \emph{elliptic shifted factorials} by
\begin{subequations}\label{esf}
\begin{equation}\label{esfa}(a;q)_k:=\theta(a)\theta(aq)\dotsm\theta(aq^{k-1}),
\end{equation}
and write
\begin{equation}(a_1,\dots,a_n;q)_k:=(a_1;q)_k\dotsm(a_n;q)_k.\end{equation}
\end{subequations}
These symbols satisfy similar identities as in the case $p=0$
\cite[Appendix I]{bhs}. In particular, we mention that
\begin{equation}\label{epdi}\frac{(a;q)_{n-k}}{(b;q)_{n-k}}=\left(\frac
    ba\right) ^k\frac{(a;q)_n(q^{1-n}/b;q)_k}{(b;q)_n(q^{1-n}/a;q)_k},
\end{equation}
and
\begin{equation}\label{epi}\frac{(a;q)_{n}}{(b;q)_{n}}=\left(\frac
    ab\right) ^n\frac{(q^{1-n}/a;q)_n}{(q^{1-n}/b;q)_n}.
\end{equation}

\section{Warnaar's matrix inversion and elliptic partial fractions}
\label{sec:pf}

In this section we give an easy proof of Warnaar's matrix inversion.
Since the case $n=l$ is trivial, it is enough to prove that
 the left-hand side of \eqref{wmi} vanishes for $n> l$.
Writing this side out explicitly gives
\begin{multline*}
\sum_{k=l}^n\frac{\prod_{j=k}^{n-1}\ta(a_jc_k)\ta(a_j/c_k)}
{\prod_{j=k+1}^{n}\ta(c_jc_k)\ta(c_j/c_k)}
\frac{c_l\,\ta(a_lc_l)\ta(a_l/c_l)}{c_k\,\ta(a_kc_k)\ta(a_k/c_k)}
\frac{\prod_{j=l+1}^{k}\ta(a_jc_k)\ta(a_j/c_k)}
{\prod_{j=l}^{k-1}\ta(c_jc_k)\ta(c_j/c_k)}\\
=c_l\theta(a_lc_l)\ta(a_l/c_l)
\sum_{k=l}^n\frac{1}{c_k}\frac{\prod_{j=l+1}^n\theta(a_jc_k)\ta(a_j/c_k)}
{\prod_{j=l,\,j\neq k}^n\theta(c_jc_k)\ta(c_j/c_k)}.
\end{multline*}
Thus, it is enough to prove that
\begin{equation}\label{mipf}
\sum_{k=l}^n\frac{1}{c_k}\frac{\prod_{j=l+1}^{n-1}\theta(a_jc_k)\ta(a_j/c_k)}
{\prod_{j=l,\,j\neq k}^n\theta(c_jc_k)\ta(c_j/c_k)}=0,\qquad n>l, 
\end{equation}
where (as a matter of relabeling) we may  assume $l=1$.

We are now reduced to a theta function identity of Gustafson 
\cite[Lemma 4.14]{gu}, which we write as
\begin{equation}\label{dpf}
\sum_{k=1}^n\frac{a_k\prod_{j=1}^{n-2}\theta(a_kb_j)\ta(a_k/b_j)}
{\prod_{j=1,\,j\neq k}^n\theta(a_ka_j)\ta(a_k/a_j)}=0, \qquad n\geq
2.\end{equation}
The case $p=0$ is equivalent to an
 elementary partial fraction expansion, so we refer to \eqref{dpf}
as an elliptic partial fraction identity.
To identify \eqref{mipf} with \eqref{dpf} it is enough to
replace $c_j$ with $a_j$, $a_j$ with $b_{j-1}$ and then use
\eqref{tif} repeatedly.

Gustafson's proof of \eqref{dpf} uses Liouville's theorem and is thus 
analytic in nature. We refer to \cite{rosers} for an elementary proof 
(using only \eqref{tif} and \eqref{tfaddf}), as well as some further
comments on this identity.

\section{An operator proof of Warnaar's matrix inversion}
\label{sec:op}

In \cite{krattopm} Krattenthaler gave a method for solving Lagrange
inversion problems, which are closely connected with the problem of inverting
lower-tri\-angular matrices. In particular, Krattenthaler applied this method
in \cite{krattmi} to derive a very general matrix inversion, namely, 
the $p=0$ case of \eqref{fg}.
In the following, we provide a proof of Warnaar's elliptic matrix inversion
using Krattenthaler's operator method. Like in Warnaar's proof, the  
essential ingredient is the addition formula \eqref{tfaddf}.

By a {\em formal Laurent series} we mean a series of the form
$\sum_{n\ge k} a_n z^{n}$, for some $k\in\mathbb Z$.
Given the formal Laurent series $a(z)$ and $b(z)$, we introduce
the bilinear form $\langle\, ,\, \rangle$ by
\begin{equation*}
\langle a(z),b(z)\rangle = \big[z^0\big](a(z)\cdot b(z)),
\end{equation*}
where $\big[z^0\big] c(z)$ denotes the coefficient
of $z^0$ in $c(z)$. Given any linear operator
$L$ acting on formal Laurent series, $L^*$ denotes the adjoint of $L$ with
respect to $\langle\,,\,\rangle$; i.e.
$\langle L a(z),b(z)\rangle=\langle a(z),L^* b(z)\rangle$
for all formal Laurent series
$a(z)$ and $b(z)$. We need the following special case of
\cite[Theorem~ 1]{krattopm}.

\begin{Lemma}\label{lemma}
Let $F=(f_{n k})_{n,k\in\mathbb Z}$ be an infinite
lower-triangular matrix with $f_{k k}\ne 0$ for all
$k\in\mathbb Z$. For $k\in\mathbb Z$, define the formal Laurent
series $f_{k}(z)$ and $g_{k}(z)$ by $f_{k}(z)=\sum_{n\ge k} f_{n k}z^n$
and $g_{k}(z)=\sum_{l\le k}g_{k l}z^{-l}$, where $(g_{k l})_{k,l\in\mathbb Z}$
is the uniquely determined inverse matrix of $F$. Suppose that
for $k\in\mathbb Z$ a system of equations of the form
\begin{equation}\label{lemmid}
U f_{k}(z)=w_k V f_{k}(z)
\end{equation}
holds, where $U$, $V$ are linear operators acting on formal Laurent series,
$V$ being bijective, and $(w_k)_{k\in\mathbb Z}$ is an
arbitrary sequence of different non-zero constants.
Then, if $h_{k}(z)$ is a solution of the dual system
\begin{equation}\label{lemmid1}
U^* h_{k}(z)=w_k V^* h_{k}(z),
\end{equation}
with $h_{k}(z)\not\equiv 0$ for all $k\in\mathbb Z$,
the series $g_{k}(z)$ is given by
\begin{equation}\label{lemmid2}
g_{k}(z)=\frac{1}{\langle f_{k}(z),
V^*h_{k}(z)\rangle} V^* h_{k}(z).
\end{equation}
\end{Lemma}

In order to prove Warnaar's elliptic extension of Krattenthaler's
matrix inversion \eqref{fg}, we set $f_k(z)=\sum_{n\ge k}f_{nk}z^k$
with $f_{nk}$ given as in \eqref{fnk}. Obviously, for $n\ge k$,
\begin{equation}\label{id1}
\ta(c_nc_k,c_n/c_k)f_{nk}=
\ta(a_{n-1}c_k,a_{n-1}/c_k)f_{n-1,k}.
\end{equation}
We now introduce a ``multiplier'' after which we apply the addition formula
for theta functions and separate the variables depending on $n$ and on $k$
appearing in \eqref{id1}. Namely, we multiply both sides of
\eqref{id1} by $\ta(uv,u/v)$ where $u,v$ are two new auxiliary
independent variables, which gives
\begin{equation}\label{id2}
\ta(c_nc_k,c_n/c_k,uv,u/v)f_{nk}=
\ta(a_{n-1}c_k,a_{n-1}/c_k,uv,u/v)f_{n-1,k}.
\end{equation}
Next, we apply the addition formula \eqref{tfaddf} to each side of
\eqref{id2} and obtain
\begin{multline*}
\Big[\ta(c_nv,c_n/v,uc_k,u/c_k)
+\frac u{c_k}\,\ta(vc_k,c_k/v,c_nu,c_n/u)\Big]f_{nk}\\=
\Big[\ta(a_{n-1}v,a_{n-1}/v,uc_k,u/c_k)
+\frac u{c_k}\,\ta(vc_k,c_k/v,a_{n-1}u,a_{n-1}/u)\Big]f_{n-1,k}.
\end{multline*}
If we define the linear operators $\mathcal A$ and $\mathcal C$ by
$\mathcal A z^k=a_kz^k$ and $\mathcal C z^k=c_kz^k$, for all $k\in\mathbb Z$,
this may be rewritten in the form
\begin{multline*}
\Big[\ta(\mathcal Cv,\mathcal C/v,uc_k,u/c_k)
+\frac u{c_k}\,\ta(vc_k,c_k/v,\mathcal Cu,\mathcal C/u)\Big]f_{k}(z)\\=
z\Big[\ta(\mathcal Av,\mathcal A/v,uc_k,u/c_k)
+\frac u{c_k}\,\ta(vc_k,c_k/v,\mathcal Au,\mathcal A/u)\Big]f_{k}(z),
\end{multline*}
or, equivalently,
\begin{multline}\label{id3}
\big[\ta(\mathcal Cv,\mathcal C/v)
-z\,\ta(\mathcal Av,\mathcal A/v)\big]f_{k}(z)\\=
\frac{u\,\ta(vc_k,c_k/v)}
{c_k\,\ta(uc_k,u/c_k)}
\big[z\,\ta(\mathcal Au,\mathcal A/u)-
\ta(\mathcal Cu,\mathcal C/u)
\big]f_{k}(z),
\end{multline}
valid for all $k\in\mathbb Z$. 

Equation \eqref{id3} is a system of equations
of type \eqref{lemmid} with
\begin{equation*}
U=\ta(\mathcal Cv,\mathcal C/v)
-z\,\ta(\mathcal Av,\mathcal A/v),
\end{equation*}
\begin{equation*}
V=z\,\ta(\mathcal Au,\mathcal A/u)-
\ta(\mathcal Cu,\mathcal C/u),
\end{equation*}
and
\begin{equation*}
w_k=\frac{u\,\ta(vc_k,c_k/v)}
{c_k\,\ta(uc_k,u/c_k)}.
\end{equation*}
The dual equations \eqref{lemmid1} for the auxiliary formal Laurent series
$h_k(z)=\sum_{l\leq k}h_{kl}z^{-l}$ in this case read
\begin{multline}\label{id4}
\big[\ta(\mathcal C^*v,\mathcal C^*/v)
-\ta(\mathcal A^*v,\mathcal A^*/v)\,z\big]h_{k}(z)\\=
\frac{u\,\ta(vc_k,c_k/v)}
{c_k\,\ta(uc_k,u/c_k)}
\big[\ta(\mathcal A^*u,\mathcal A^*/u)\,z-
\ta(\mathcal C^*u,\mathcal C^*/u)
\big]h_{k}(z).
\end{multline}
Since $\mathcal A^*z^{-k}=a_kz^{-k}$ and
$\mathcal C^*z^{-k}=c_kz^{-k}$, by comparing coefficients of
$z^{-l}$ in \eqref{id4} we obtain
\begin{multline*}
\Big[\ta(c_lv,c_l/v,uc_k,u/c_k)
+\frac u{c_k}\,\ta(vc_k,c_k/v,c_lu,c_l/u)\Big]h_{kl}\\=
\Big[\ta(a_{l}v,a_{l}/v,uc_k,u/c_k)
+\frac u{c_k}\,\ta(vc_k,c_k/v,a_{l}u,a_{l}/u)
\Big]h_{k,l+1},
\end{multline*}
which, after application of the addition formula \eqref{tfaddf}
and dividing both sides by $\ta(uv,u/v)$,
reduces to
\begin{equation*}
\ta(c_lc_k,c_l/c_k)h_{kl}=
\ta(a_{l}c_k,a_{l}/c_k)h_{k,l+1}.
\end{equation*}
If we set $h_{kk}=1$, we get
\begin{equation*}
h_{kl}=
\frac{\prod_{j=l}^{k-1}\ta(a_jc_k,a_j/c_k)}
{\prod_{j=l}^{k-1}\ta(c_jc_k,c_j/c_k)}.
\end{equation*}

Taking into account \eqref{lemmid2}, we compute
\begin{multline}\label{id5}
V^*h_k(z)=\big[\ta(\mathcal A^*u,\mathcal A^*/u)\,z-
\ta(\mathcal C^*u,\mathcal C^*/u)\big]h_k(z)\\
=\sum_{l\le k}
\bigg[\frac{\ta(c_lc_k,c_l/c_k)}{\ta(a_lc_k,a_l/c_k)}
\ta(a_lu,a_l/u)-\ta(c_lu,c_l/u)\bigg]
\frac{\prod_{j=l}^{k-1}\ta(a_jc_k,a_j/c_k)}
{\prod_{j=l}^{k-1}\ta(c_jc_k,c_j/c_k)}\,z^{-l}\\
=\sum_{l\le k}\ta(c_kv,c_k/v)
\frac{a_l\,\ta(a_lc_l,c_l/a_l)}
{c_k\,\ta(a_lc_k,a_l/c_k)}\,
\frac{\prod_{j=l}^{k-1}\ta(a_jc_k,a_j/c_k)}
{\prod_{j=l}^{k-1}\ta(c_jc_k,c_j/c_k)}\,z^{-l},
\end{multline}
where we again have used the addition formula \eqref{tfaddf}.
Now, since $f_{kk}=1$, the pairing $\langle f_{k}(z),V^*h_{k}(z)\rangle$
is simply the coefficient of $z^{-k}$ in \eqref{id5}.
Thus, \eqref{lemmid2} reads
\begin{equation}\label{id6}
g_k(z)=-\frac{1}{\ta(c_kv,c_k/v)}V^*h_k(z)
\end{equation}
where $g_k(z)=\sum_{l\le k}g_{kl}z^{-l}$. Hence, extracting coefficients of
$z^{-l}$ in \eqref{id6} we obtain exactly \eqref{gkl}.

\section{Elliptic Karlsson--Minton-type identities}
\label{sec:km}

As was mentioned in the introduction, we can obtain a 
generalization of the Karlsson--Minton-type identities \eqref{wbb} and
\eqref{kmsi} as a special case of the partial fraction identity
\eqref{dpf}. To this end, we make the substitutions
\begin{subequations}\label{sub}
\begin{equation}
(a_1,\dots,a_n)\mapsto(a_1,a_1q,\dots,a_1q^{l_1},\dots,a_s,a_sq,\dots,
a_sq^{l_s}),
\end{equation}
\begin{multline}(b_1,\dots,b_{n-2})\\
\mapsto(b_1,b_1q^{1/y_1},\dots,b_1q^{(m_1-1)/y_1},\dots,
b_r,b_rq^{1/y_r},\dots,b_rq^{(m_r-1)/y_r})\end{multline}
\end{subequations}
in \eqref{dpf}, with $m_i$ and $l_i$ non-negative and $y_i$ positive integers
satisfying 
\begin{equation}\label{parcond}|l|+s=|m|+2.\end{equation}
The resulting special case of   \eqref{dpf} may be written
$$\sum_{i=1}^s\sum_{k=0}^{l_i}\frac{a_iq^k\prod_{j=1}^r
\prod_{t=0}^{m_j-1}\theta(a_iq^kb_jq^{t/y_j},a_iq^k/b_jq^{t/y_j})}
{\prod_{t=0,\,t\neq
      k}^{l_i}\theta(a_i^2q^{k+t},q^{k-t})\prod_{j=1,\,j\neq
    i}^s\prod_{t=0}^{l_j}\theta(a_iq^ka_jq^t,a_iq^k/a_jq^t) }=0.$$
It is now straight-forward to rewrite the products in $t$ in terms of
elliptic shifted factorials, giving
$$\prod_{t=0}^{m_j-1}\theta(a_iq^kb_jq^{t/y_j})=
(a_ib_jq^k;q^{1/y_j})_{m_j}=
(a_ib_j;q^{1/y_j})_{m_j}\frac{(a_ib_jq^{m_j/y_j};q^{1/y_j})_{y_jk}}
{(a_ib_j;q^{1/y_j})_{y_jk}}, $$
and similarly
$$\prod_{t=0}^{m_j-1}\theta(a_iq^k/b_jq^{t/y_j})
=(a_iq^{(1-m_j)/y_j}/b_j;q^{1/y_j})_{m_j}
\frac{(a_iq^{1/y_j}/b_j;q^{1/y_j})_{y_jk}}
{(a_iq^{(1-m_j)/y_j}/b_j;q^{1/y_j})_{y_jk}},$$
$$\frac 1{\prod_{t=0,\,t\neq
      k}^{l_i}\theta(a_i^2q^{k+t})}=
\frac{1}{(a_i^2q;q)_{l_i}}\frac{\theta(a_i^2q^{2k})}{\theta(a_i^2)}
\frac{(a_i^2;q)_k}{(a_i^2q^{l_i+1};q)_k},$$
$$\frac1{\prod_{t=0,\,t\neq
      k}^{l_i}\theta(q^{k-t})}=\frac 1{(q^{-l_i};q)_{l_i}}
\frac{(q^{-l_i};q)_k}{(q;q)_k},$$
$$\frac 1{\prod_{t=0}^{l_j}\theta(a_iq^ka_jq^t)}=
\frac{1}{(a_ia_j;q)_{l_j+1}}\frac{(a_ia_j;q)_k}{(a_ia_jq^{l_j+1};q)_k},$$
$$\frac 1{\prod_{t=0}^{l_j}\theta(a_iq^k/a_jq^t)}
=\frac{1}{(a_iq^{-l_j}/a_j;q)_{l_j+1}}\frac{(a_iq^{-l_j}/a_j;q)_k}
{(a_iq/a_j;q)_k}.$$

We thus arrive at the following result.

\begin{Theorem}\label{kmt}
Let $l_1,\dots,l_s$  and $m_1,\dots,m_r$ be non-negative integers such
that $|l|+s=|m|+2$, and let $y_1,\dots,y_r$ be positive integers. Then
the following identity holds:
\begin{multline*}\sum_{i=1}^s \frac{a_i\prod_{j=1}^r(a_ib_j,
a_iq^{(1-m_j)/y_j}/b_j;q^{1/y_j})_{m_j}}{(a_i^2q,q^{-l_i};q)_{l_i}
\prod_{j=1,\,j\neq i}^s(a_ia_j,a_iq^{-l_j}/a_j;q)_{l_j+1}}\\
\times\sum_{k=0}^{l_i} \frac{\theta(a_i^2q^{2k})}{\theta(a_i^2)}\,q^k
\prod_{j=1}^s\frac{(a_ia_j,a_iq^{-l_j}/a_j;q)_k}
{(a_iq/a_j,a_ia_jq^{l_j+1};q)_k}\\
\times\prod_{j=1}^r\frac
{(a_ib_jq^{m_j/y_j}, a_iq^{1/y_j}/b_j;q^{1/y_j})_{y_jk}}
{( a_iq^{(1-m_j)/y_j}/b_j,a_ib_j;q^{1/y_j})_{y_jk}}=0.
\end{multline*}
\end{Theorem}

\begin{Remark}
It is clear from the proof that Theorem \ref{kmt} is actually equivalent
to its special case when $y_j\equiv 1$. This may be checked directly
using \eqref{xsk}. However, in view of the work of Warnaar
\cite{warnell}, the form given above 
seems more useful for potential application to quadratic and
higher identities. 
\end{Remark}

\begin{Remark}
In principle, one can obtain an even more general identity by
replacing \eqref{sub} with a substitution involving independent bases,
that is,
\begin{gather*}
(a_1,\dots,a_n)\mapsto(a_1,\dots,a_1q_1^{l_1},\dots,a_s,\dots,
a_sq_s^{l_s}),\\
(b_1,\dots,b_{n-2})\mapsto(b_1,\dots,b_1p_1^{m_1-1},\dots,
b_r,\dots,b_rp_r^{m_r-1}).\end{gather*}
However, the inner sums in the resulting identity will not be elliptic 
hypergeometric.
\end{Remark}

\begin{Remark}\label{searsremark}
In the basic case, $p=0$, Theorem \ref{kmt} may be obtained as a special
case of Sears' transformation for well-poised series
\cite{sears}. More precisely, if we start from the special case given
  in \cite[Exercise 4.7]{bhs}, replace $r$ by $r+s$ and choose the
parameters $(b_1,\dots,b_{r+s})$ there as
$$(q^{-l_1}/a_1,\dots,q^{-l_s}/a_s,q^{m_1+1}/a_{s+1},\dots,q^{m_r+1}/a_{r+s}),
$$
we obtain an identity equivalent to the case $p=0$ of
 Theorem \ref{kmt}. This is exactly
the case of Sears' transformation 
 when all series involved  are
terminating, very-well-poised and balanced. Since these restrictions
are natural in the elliptic case \cite{spir}, we may view Theorem
\ref{kmt} as an elliptic analogue of Sears' transformation. 
\end{Remark}

For applications, the case $s=2$ of Theorem \ref{kmt} seems especially
useful, and we give it explicitly in the following corollary. 
We have made the substitutions
$(a_1,a_2,l_1,l_2,b_j)\mapsto(\sqrt a,b/\sqrt a,N,L,c_j/\sqrt a)$
and used \eqref{epi} to simplify some of the factors.

\begin{Corollary}\label{trc}
Let $L$, $N$ and $m_1,\dots,m_r$ be non-negative integers with 
$|m|=L+N$, and let $y_1,\dots,y_r$ be positive integers. Then,
\begin{multline*}
\sum_{k=0}^N\frac{\theta(aq^{2k})}{\theta(a)}
\frac{(a,q^{-N},b,aq^{-L}/b;q)_k}{(q,aq^{N+1},aq/b,bq^{L+1};q)_k}\,q^k
\prod_{j=1}^r\frac{(c_jq^{m_j/y_j},aq^{1/y_j}/c_j;q^{1/y_j})_{y_jk}}
{(aq^{(1-m_j)/y_j}/c_j,c_j;q^{1/y_j})_{y_jk}}\\
=\frac{(aq,q;q)_N}{(bq,aq/b;q)_N}\frac{(bq,bq/a;q)_L}{(b^2q/a,q;q)_L}
\prod_{j=1}^r
\frac{(c_j/b,c_jb/a;q^{1/y_j})_{m_j}}
{(c_j,c_j/a;q^{1/y_j})_{m_j}}\\
\times\sum_{k=0}^L\frac{\theta(b^2q^{2k}/a)}{\theta(b^2/a)}
\frac{(b^2/a,q^{-L},b,bq^{-N}/a;q)_k}{(q,q^{L+1}b^2/a,bq/a,bq^{N+1};q)_k}\,q^k
\\
\times\prod_{j=1}^r\frac{(bc_jq^{m_j/y_j}/a,bq^{1/y_j}/c_j;q^{1/y_j})_{y_jk}}
{(bq^{(1-m_j)/y_j}/c_j,bc_j/a;q^{1/y_j})_{y_jk}}.
\end{multline*}
\end{Corollary}

If we let $L=0$ in  Corollary \ref{trc} we obtain the following
summation formula.

\begin{Corollary}\label{mbkms}
Let $y_1,\dots,y_r$ be positive integers and $m_1,\dots,m_r$ be
non-negative integers with $m_1+\dots+m_r=N$.
Then the following identity holds:
\begin{multline*}\sum_{k=0}^N\frac{\theta(aq^{2k})}{\theta(a)}
\frac{(a,q^{-N},b,a/b;q)_k}{(q,aq^{N+1},aq/b,bq;q)_k}\,q^k
\prod_{j=1}^r\frac{(c_jq^{m_j/y_j},aq^{1/y_j}/c_j;q^{1/y_j})_{y_jk}}
{(aq^{(1-m_j)/y_j}/c_j,c_j;q^{1/y_j})_{y_jk}}\\
=\frac{(aq,q;q)_N}{(bq,aq/b;q)_N}\prod_{j=1}^r
\frac{(c_j/b,c_jb/a;q^{1/y_j})_{m_j}}
{(c_j,c_j/a;q^{1/y_j})_{m_j}}. \end{multline*}
\end{Corollary}

Note that the case $r=1$ of Corollary \ref{mbkms} is equivalent
to \eqref{wbb}, and that the case
$y_j\equiv 1$ is   \eqref{kmsi}.

\section{Some exotic Karlsson--Minton-type identities}
\label{sec:ekm}

Besides \eqref{dpf}, 
we are aware of another
 elliptic partial fraction expansion, namely,
\begin{equation}\label{apf}\sum_{k=1}^n\frac{\prod_{j=1}^n\theta(a_k/b_j)}
{\prod_{j=1,\,j\neq k}^n\theta(a_k/a_j)}=0,\qquad a_1\dotsm a_n=b_1\dotsm b_n,
\end{equation}
 which goes back at least to the 1898 treatise of Tannery and  Molk
\cite[p.~34]{tm}. Again, we refer to \cite{rosers} for an elementary
proof and some further comments. 

It does not seem possible to obtain a
matrix inversion from \eqref{apf} in a similar way as Warnaar's 
inversion was obtained from \eqref{dpf} in Section
\ref{sec:pf}. However, it is straight-forward to imitate the analysis
of Section \ref{sec:km} and obtain Karlsson--Minton-type summation
and transformation formulas from \eqref{apf}. The resulting identities seem
quite exotic and appear to be new even in the case $p=0$.

Thus, we make the substitutions \eqref{sub} into \eqref{apf}. In place
of \eqref{parcond}, we now have the two conditions $|l|+s=|m|$ and
\begin{equation}\label{apc}q^{\binom{l_1+1}{2}+\dots+\binom{l_s+1}{2}}
a_1^{l_1+1}\dotsm
a_s^{l_s+1}
=q^{\frac{1}{y_1}\binom{m_1}{2}+\dots
+\frac{1}{y_r}\binom{m_r}{2}}b_1^{m_1}\dotsm b_r^{m_r}.\end{equation}
Clearly, the resulting transformation can  be obtained from
Theorem \ref{kmt} by deleting the factor $a_iq^k$, together with all
factors involving products (rather than quotients) of the parameters
$a_i$, $b_i$. This gives the following result.

\begin{Theorem}\label{atr}
Let $l_1,\dots,l_s$  and $m_1,\dots,m_r$ be non-negative integers 
and  $y_1,\dots,y_r$ be positive integers. Assume that
 $|l|+s=|m|$, and that \eqref{apc} holds. Then, 
\begin{multline*}\sum_{i=1}^s \frac{\prod_{j=1}^r
(a_iq^{(1-m_j)/y_j}/b_j;q^{1/y_j})_{m_j}}{(q^{-l_i};q)_{l_i}
\prod_{j=1,\,j\neq i}^s(a_iq^{-l_j}/a_j;q)_{l_j+1}}\\
\times\sum_{k=0}^{l_i}
\prod_{j=1}^s\frac{(a_iq^{-l_j}/a_j;q)_k}
{(a_iq/a_j;q)_k}\prod_{j=1}^r\frac
{(a_iq^{1/y_j}/b_j;q^{1/y_j})_{y_jk}}
{(a_iq^{(1-m_j)/y_j}/b_j;q^{1/y_j})_{y_jk}}=0.
\end{multline*}
\end{Theorem}

Next we write down the case $s=2$ of  Theorem \ref{atr} explicitly. 
 For this we  make the substitutions 
$$(a_1,a_2,l_1,l_2,b_j)\mapsto(bq^L,1,N,L,
bq^{(Ly_j-m_j+1)/y_j}/c_j).$$
 (Since we may multiply all $a_j$ and $b_j$ in
\eqref{apf} with a common factor, the assumption
$a_2=1$ is no restriction.) This yields that if $|m|=N+L+2$ and
\begin{equation}\label{npc}
q^{\binom{L+1}2}b^{L+1}=q^{\binom{N+1}2+\frac{1}{y_1}\binom{m_1}{2}+\dots
+\frac{1}{y_r}\binom{m_r}{2}}c_1^{m_1}\dotsm c_r^{m_r},\end{equation} 
then
\begin{multline*}\frac{\prod_{j=1}^r(c_j;q^{1/y_j})_{m_j}}
{(q^{-N};q)_N(b;q)_{L+1}}
\sum_{k=0}^N
\frac{(q^{-N},b;q)_k}{(q,bq^{L+1};q)_k}
\prod_{j=1}^r\frac{(c_jq^{m_j/y_j};q^{1/y_j})_{y_jk}}{(c_j;q^{1/y_j})_{y_jk}}
\\
+\frac{\prod_{j=1}^r(q^{-L}c_j/b;q^{1/y_j})_{m_j}}
{(q^{-L};q)_L(q^{-L-N}/b;q)_{N+1}}\sum_{k=0}^L
\frac{(q^{-L},q^{-L-N}/b;q)_k}{(q,q^{1-L}/b;q)_k}\\
\times\prod_{j=1}^r\frac{(q^{(m_j-Ly_j)/y_j}c_j/b;q^{1/y_j})_{y_jk}}
{(q^{-L}c_j/b;q^{1/y_j})_{y_jk}}=0.
\end{multline*}
To make this look nicer  we replace $k$ by $L-k$ in the second sum. After
repeated application of \eqref{epdi} and some further simplification,
we arrive at the following  transformation
formula.

\begin{Corollary}\label{akmt}
Let $L$, $N$ and $m_1,\dots,m_r$ be non-negative integers
with $|m|=N+L+2$,  
and let  $y_1,\dots,y_r$ be positive integers. Then, assuming also \eqref{npc},
one has the identity
\begin{multline*}
\sum_{k=0}^N
\frac{(q^{-N},b;q)_k}{(q,bq^{L+1};q)_k}
\prod_{j=1}^r\frac{(c_jq^{m_j/y_j};q^{1/y_j})_{y_jk}}
{(c_j;q^{1/y_j})_{y_jk}}\\
=b^{N+1}\frac{(q;q)_N(bq;q)_L}{(bq;q)_N(q;q)_L}\prod_{j=1}^r
\frac{(c_j/b;q^{1/y_j})_{m_j}}
{(c_j;q^{1/y_j})_{m_j}}\\
\times\sum_{k=0}^L
\frac{(q^{-L},b;q)_k}{(q,bq^{N+1};q)_k}
\prod_{j=1}^r\frac{(bq^{1/y_j}/c_j;q^{1/y_j})_{y_jk}}
{(bq^{(1-m_j)/y_j}/c_j;q^{1/y_j})_{y_jk}}.
\end{multline*}
\end{Corollary}

When $L=0$, Corollary \ref{akmt} reduces to the following summation formula.

\begin{Corollary}\label{akms}
Let  $m_1,\dots,m_r$ and $N$ be non-negative integers with $|m|=N+2$,
and let  $y_1,\dots,y_r$  be positive integers. Then, assuming also 
$$b=q^{\binom{N+1}{2}+
\frac{1}{y_1}\binom{m_1}{2}+\dots
+\frac{1}{y_r}\binom{m_r}{2}}c_1^{m_1}\dotsm c_r^{m_r}, $$
one has the identity
$$
\sum_{k=0}^N
\frac{(q^{-N},b;q)_k}{(q,bq;q)_k}
\prod_{j=1}^r\frac{(c_jq^{m_j/y_j};q^{1/y_j})_{y_jk}}
{(c_j;q^{1/y_j})_{y_jk}}=b^{N+1}\frac{(q;q)_N}{(bq;q)_N}\prod_{j=1}^r
\frac{(c_j/b;q^{1/y_j})_{m_j}}
{(c_j;q^{1/y_j})_{m_j}}.
$$
\end{Corollary}

The evaluation in Corollary \ref{akms}
looks so unusual that it is worth pointing out that we
 believe it is free from misprints. In particular, the factor
$q^k$ is not missing from the left-hand side. The case
 $p=0$, which also appears  to be new,  
has been confirmed by numerical calculations. 

To our knowledge, Corollary \ref{akms} is the first known summation
formula for a theta hypergeometric series (in the terminology of 
Spiridonov \cite{spir}) that is not well-poised, and also the first
such result with argument $\neq q$. On the other hand, using the
results of  \cite{spir} it is easy to 
 check that the sum  is \emph{modular} (that is, 
  invariant under a
 natural action of $\mathrm{SL}(2,\mathbb Z)$  on $(p,q)$-space)
and, 
in particular,  \emph{balanced} in the sense of
Spiridonov.  However, the special case $p=0$ is not balanced in the
 usual sense of basic hypergeometric series  \cite{bhs}. 
This is another  indication of the 
 importance of modular invariance and Spiridonov's balancing condition
for elliptic hypergeometric series.

\section*{Appendix. An alternative proof of \eqref{kmsi}}
\setcounter{section}{1}
\renewcommand{\thesection}{\Alph{section}}
\setcounter{equation}{0}

 Gasper's proof of the case $p=0$
of \eqref{kmsi} uses induction on $r$. As was remarked in the introduction,
this proof does not immediately extend to the general case. However,
we were able to find a proof by induction on $N$, which is different
in  details from Gasper's proof, but closer to standard methods for
basic hypergeometric series  \cite{bhs}
than the proof given in Section~\ref{sec:km}. We include this proof
here both since it may have independent interest and since it may be
useful for generalizations, for instance, to multiple series. For
brevity, we will write $(a)_k=(a;q)_k$ for the elliptic shifted factorials
in \eqref{esfa}.

To start our inductive proof of \eqref{kmsi}, we assume that it
 holds for
fixed $N$ and consider the sum
$$S=\sum_{k=0}^{N+1}\frac{\theta(aq^{2k})}{\theta(a)}
\frac{(a,q^{-N-1},b,a/b)_k}{(q,aq^{N+2},aq/b,bq)_k}\,q^k
\prod_{j=1}^r\frac{(c_jq^{m_j},aq/c_j)_k}{(aq^{1-m_j}/c_j,c_j)_k},$$
where $m_1+\dots+m_r=N+1$. By symmetry, we may assume $m_r\geq 1$.

We multiply the sum $S$ termwise by
\begin{multline*}1=
\frac{1}
{\ta(aq^{N+1},q^{-N-1},c_rq^{m_r+k-1},aq^{1-m_r+k}/c_r)}\\
\times\big[\ta(aq^{k+N+1},q^{k-N-1},c_rq^{m_r-1},aq^{1-m_r}/c_r)\\
-q^{-N-1}
\ta(aq^{k},q^{k},c_rq^{m_r+N},aq^{2-m_r+N}/c_r) \big],
\end{multline*}
which is equivalent to \eqref{tfaddf} with the replacements 
$$(u,v,x,y)\mapsto(\sqrt a,q^{m_r-1}c_r/\sqrt a,q^{N+1}\sqrt a,q^k\sqrt a).$$
Since the factors $\ta(q^{k-N-1})$ and $\theta(q^k)$ vanish at the
end-points $k=N+1$ and $k=0$, respectively, this gives an identity of
the form 
$$S=\sum_{k=0}^N(\dotsm)+\sum_{k=1}^{N+1}(\dotsm).$$
Replacing $k$ by $k+1$ in the last sum and simplifying gives
\begin{multline*}S=\sum_{k=0}^N\frac{\theta(aq^{2k})}{\theta(a)}
\frac{(a,q^{-N},b,a/b,c_rq^{m_r-1},aq/c_r)_k}
{(q,aq^{N+1},aq/b,bq,aq^{2-m_r}/c_r,c_r)_k}\,q^k
\prod_{j=1}^{r-1}\frac{(c_jq^{m_j},aq/c_j)_k}
{(aq^{1-m_j}/c_j,c_j)_k}\\ 
-q^{-N}\frac{\ta(aq,aq^2,b,a/b,c_rq^{m_r+N},aq/c_r,aq^{2-m_r+N}/c_r)}
{\ta(aq/b,bq,aq^{N+1},aq^{N+2},aq^{1-m_r}/c_r,aq^{2-m_r}/c_r,c_r)}
\prod_{j=1}^{r-1}\frac{\ta(c_jq^{m_j},aq/c_j)}{\ta(aq^{1-m_j}/c_j,c_j)}
\\
\times\sum_{k=0}^N\frac{\theta(aq^{2k+2})}{\theta(aq^2)}
\frac{(aq^2,q^{-N},bq,aq/b,c_rq^{m_r},aq^{2}/c_r)_k}
{(q,aq^{N+3},aq^2/b,bq^2,aq^{3-m_r}/c_r,c_rq)_k}\,q^k\\
\times
\prod_{j=1}^{r-1}\frac{(c_jq^{1+m_j},aq^{2}/c_j)_k}{(aq^{2-m_j}/c_j,c_jq)_k}.
\end{multline*}

Both sums are now evaluated by the induction hypothesis, giving
\begin{multline*}
S=\frac{(aq,q)_N}{(bq,aq/b)_N}\frac{(c_r/b,c_rb/a)_{m_r-1}}
{(c_r,c_r/a)_{m_r-1}}
\prod_{j=1}^{r-1}\frac{(c_j/b,c_jb/a)_{m_j}}
{(c_j,c_j/a)_{m_j}}\\
\times\left\{1-\frac{\ta(b,a/b,c_rq^{N+m_r},aq^{N+2-m_r}/c_r)}
{\ta(bq^{N+1},aq^{N+1}/b,c_rq^{m_r-1},aq^{1-m_r}/c_r)}\right\}.
\end{multline*}
Using again \eqref{tfaddf}, this time with 
$$(u,v,x,y)\mapsto(\sqrt a,q^{m_r-1}c_r/\sqrt a,q^{N+1}\sqrt a,b/\sqrt a), $$
we find that the factor within brackets equals
$$\frac{\ta(q^{N+1},aq^{N+1},q^{m_r-1}c_r/b,q^{m_r-1}c_rb/a)}
{\ta(bq^{N+1},aq^{N+1}/b,q^{m_r-1}c_r,q^{m_r-1}c_r/a)}, $$
and thus
$$S=\frac{(aq,q)_{N+1}}{(bq,aq/b)_{N+1}}
\prod_{j=1}^r\frac{(c_j/b,c_jb/a)_{m_j}}
{(c_j,c_j/a)_{m_j}}. $$
This completes our alternative proof of \eqref{kmsi}.

\end{document}